\documentclass[11pt]{article}

\usepackage{amsmath,amssymb,amsthm,bm}
\usepackage{bbm}
\usepackage{enumerate}
\usepackage{tikz-cd} 
\usepackage[utf8]{inputenc}
\usepackage{xcolor}
\usepackage{tikz}
\usetikzlibrary{shadings}
\usepackage{url}

\DeclareMathOperator{\sign}{sign}

\makeatletter
\providecommand*{\diff}%
{\@ifnextchar^{\DIfF}{\DIfF^{}}}
\def\DIfF^#1{%
	\mathop{\mathrm{\mathstrut d}}%
	\nolimits^{#1}\gobblespace}
\def\gobblespace{%
	\futurelet\diffarg\opspace}
\def\opspace{%
	\let\DiffSpace\!
	\ifx\diffarg(%
	\let\DiffSpace\relax
	\else
	\ifx\diffarg[%
	\let\DiffSpace\relax
	\else
	\ifx\diffarg\{%
	\let\DiffSpace\relax
	\fi\fi\fi\DiffSpace}
\makeatother

\newtheorem{theorem}{Theorem}
\newtheorem*{theorem*}{Theorem}

\newtheorem*{corollary*}{Corollary}
\newtheorem{lemma}{Lemma}
\newtheorem*{proposition*}{Proposition}

\theoremstyle{definition}

\theoremstyle{remark}

\newcommand{\ord}{k}

\DeclareMathOperator{\sgn}{sgn}

\title{On a class of sharp Sobolev type estimates with weights}
\author{Raul Hindov \and Evgeniy Lokharu}

\begin{document}
	
	\maketitle
	
	\begin{abstract}
		We study sharp weighted Sobolev-type inequalities of the form
		\[
		\int_{0}^{1}|u(x)|\rho(x) \diff x
		\leqslant \Lambda
		\Bigl(\int_{0}^{1}|u^{(k)}(x)|^2 \diff x\Bigr)^{1/2},
		\qquad u\in H_0^k(0,1),
		\]
		where $\rho$ is a non-negative weight. We characterize the minimizers and
		identify the optimal constant $\Lambda(k,\rho)$ by showing that every minimizer has a constant sign and therefore solves a nonlinear eigenvalue problem of polyharmonic type. This yields an explicit characterization of extremizers for a broad class of weights. Moreover, we even provide with a an explicit computation of the optimal constant in terms of the weight function. The new weighted estimates turn to be very useful and, among other applications, allow us to recover several previous sharp estimates and Hardy type inequalities on finite intervals.
	\end{abstract}

	\section{Introduction}
	
	Sobolev time inequalities of the form
	\[
	\|u\|_{X} \leqslant \Lambda \|u^{(k)}\|_Y
	\]
	with $X$ and $Y$ being Banach spaces of functions of one variable have been extensively studied. For a comprehensive overview of such results in one variable, we refer to \cite{Nazarov_Shcheglova2025} and \cite{Kuznetsov}. A typical example is $X = L^q(0,1)$ and $Y = L^p(0,1)$, while $u$ satisfies Dirichlet boundary conditions 
	\begin{equation}\label{eq:bc}
		u^{(j)}(0) = u^{(j)}(1) = 0, \ \ 0 \leqslant j \leqslant k-1.
	\end{equation}
	For $p=2$ the optimal constant and the extremizers are known
	in several cases:
	\begin{itemize}
		\item $k=1$ and any $q \geqslant 1$ by Schmidt \cite{Schmidt1940};
		\item any $k \geqslant 1$ and $q=2$ by Petrova \cite{Petrova2017};
		\item $k = 2,3$ and $q=\infty$ by Oshime \cite{Oshime2008} and Watanabe \cite{Watanabe2009};
		\item any $k \geqslant 1$ and $q = \infty$ by Kalyabin \cite{Kalyabin1};
		\item any $k \geqslant 1$ and $q = 1$ very recently in \cite{Hindov}.
	\end{itemize}
	Other parameter regimes for $p = 2$ are still uncovered. The purpose of the present paper is to extend the result of \cite{Hindov} to weighted spaces $X = L^1(-1,1;\rho)$ and to give a significantly simplified argument. The new weighted estimates turn to be very useful and, among other applications, allow us to recover several previous sharp estimates and some new Hardy type inequalities on finite intervals.
	
	In what follows, we will consider inequalities with weights of the form
	\begin{equation} \label{sobolev}
		\int_{0}^{1} |u(x)| \rho(x) \diff x \leqslant \Lambda \left( \int_{0}^{1} |u^{(k)}(x)|^2 \diff x  \right)^{1/2},
	\end{equation}
	where $\rho$ is a non-negative weight function, $k \geqslant 1$ is an integer, and $u \in H^k_0(0,1)$. Here, the space $H^k_0(-1,1)$ consists of all functions $u \in H^k(0,1)$ that satisfy the Dirichlet boundary conditions \eqref{eq:bc}.

	We denote by $\Lambda(k,\rho)$ the optimal constant for inequality \eqref{sobolev}, which is formally defined as
	\begin{equation} \label{inf}
		\frac{1}{\Lambda(k,\rho)} = \inf_{\substack{u \in H^k_0(0,1); \\ u \ne 0}} \frac{ \left( \int_{0}^{1} |u^{(k)}(x)|^2 \diff x  \right)^{1/2}}{\int_{0}^{1} |u(x)| \rho(x) \diff x}.
	\end{equation}
	This infimum is clearly finite and positive, and the minimizer is described in the following theorem, which is the main result of this paper.
	
	\begin{theorem}\label{thm:1}
		Let $\rho \in L^1(0,1)$ be a given non-negative function and $k \in \mathbb{N}$. Then the minimizer of \eqref{inf} corresponds to the first positive "eigenfunction" corresponding to the smallest "eigenvalue" $\mu$ of the problem
		\begin{equation} \label{eq:eigenproblem}
			\left\{
			\begin{aligned}
				(-1)^\ord u^{(2\ord)}(x) &=\mu \rho(x) \sign(u),  &&x\in [0,1], \\
				\int_{0}^{1}|u(x)| \rho(x) \diff x &= 1, \\
				u^{(j)}(0) = u^{(j)}(1) &= 0,  &&0 \leqslant j \leqslant k-1.
			\end{aligned}
			\right.
		\end{equation}		
		The corresponding unique scalar $\mu$ coincides with $\left(\frac{1}{\Lambda(k,\rho)}\right)^2$.
	\end{theorem}
	
	It is clear that the squared reciprocal of the optimal constant in \eqref{inf} corresponds to the smallest eigenvalue of the problem above. What we prove is that the corresponding eigenfunction has a sign. Then a minimizer can be found by solving a simple differential equation, where the right-hand side of the first equation of \eqref{eq:eigenproblem} is just a constant times $\rho$. The minimizer is then uniquely defined by the second normalization constraint in \eqref{eq:eigenproblem}. \\
	
	In addition to the characterization above, there is an explicit way to compute $\mu$ described as  below. Multiplying the differential equation by $u(x)$ and integrating over the entire interval yields
	\begin{equation}\label{mu}
		\mu = \int_0^1 \left( u^{(k)}(x) \right)^2 \diff x,
	\end{equation}
	by repeated integration by parts $k$ times and using all boundary conditions and the normalization constraint.
	
	After integrating the differential equation $k$ times from 0 to $x$, the $k$-th derivative has the form
	\begin{equation}\label{u_k}
		u^{(k)}(x)
		=\sum_{j=0}^{k-1} u^{(2k-1-j)}(0) \frac{x^{k-1-j}}{(k-1-j)!}
		+ (-1)^k \mu \int_0^x \frac{(x-t)^{k-1}}{(k-1)!} \, \rho(t) \diff t.
	\end{equation}
	
	We integrate the $k$-th derivative \eqref{u_k} further repeatedly $k$ times from 0 to $x$, and since we can now use conditions $u^{(j)}(0)=u^{(j)}(1)=0$,  for $0 \leqslant j \leqslant k-1$, we get a $k\times k$ linear system for coefficients $u^{(2k-1-j)}(0)$, for $0 \leqslant j \leqslant k-1$. We have
	\[
	\sum_{j=0}^{k-1}u^{(2k-1-j)}(0) \frac{(k+m)!}{(k+m-j)!}=(-1)^{k+1}\mu\int_0^1(1-t)^{k+m}\rho(t)\diff t, \quad m=0,1,\dots,k-1.
	\]
	By expressing that in matrix form, we get
	\[
	A \textbf{u}=\mu\textbf{b},
	\]
	where $A_{mj}=\frac{(k+m)!}{(k+m-j)!}$, $\textbf{u}=\left(u^{(2k-1)}(0), u^{(2k-2)}(0), \dots, u^{(k)}(0)\right)^\top\in\mathbb{R}^k$, $\textbf{b}=(b_0,b_1,\dots,b_{k-1})^\top\in\mathbb{R}^k$, and $b_m=(-1)^{k+1}\int_0^1(1-t)^{k+m}\rho(t)\diff t$.
	
	Since $A$ is a Vandermonde type matrix, it is invertible and we can solve the matrix equation for $\textbf{u}$. We have
	\[
	u^{(2k-1-j)}(0)=\mu \sum_{m=0}^{k-1}\left(A^{-1}\right)_{jm}b_m, \quad j=0,1,\dots,k-1.
	\]
	
	Using these expressions and inserting the formula \eqref{u_k} into the identity for $\mu$ in \eqref{mu}, we obtain that
	\begin{equation}
		\frac{1}{\mu}
		= \int_0^1 
		\left[
		\sum_{j=0}^{k-1}\frac{x^{k-1-j}}{(k-1-j)!} \sum_{m=0}^{k-1}\left(A^{-1}\right)_{jm}b_m 
		+ (-1)^k \int_0^x \frac{(x-t)^{k-1}}{(k-1)!} \, \rho(t) \diff t
		\right]^2 \diff x.
	\end{equation}

	Then, for example, in the case $k=1$,
	\begin{equation*}
		-u'' = \mu \rho,
		\quad \int_0^1 u \rho = 1,
		\quad u(0) = u(1) = 0,
	\end{equation*}
	we obtain
	\begin{equation*}
		\frac{1}{\mu} = \int_0^1 \left( C_0 - \int_0^x \rho(t) \, dt \right)^2 \diff x,
	\end{equation*}
	where $C_0 = \int_0^1 (1-t) \rho(t) \diff t$.
	
	In the case $k=2$,
	\begin{equation*}
		u^{(4)} = \mu \rho,
		\quad \int_0^1 u \rho = 1,
		\quad u(0)=u(1)=u'(0)=u'(1)=0,
	\end{equation*}
	we obtain
	\begin{equation*}
		\frac{1}{\mu} = \int_0^1 \left( C_0 + C_1 x + 
		\int_0^x (x-t) \rho(t) \, \diff t \right)^2 \diff x,
	\end{equation*}
	where
	\begin{align*}
		C_0 &= \int_0^1 (1-t)^{2}\,\rho(t) \diff t - \int_0^1 (1-t)^{3}\,\rho(t) \diff t \\[6pt]
		C_1 &= -3\int_0^1 (1-t)^{2}\,\rho(t) \diff t +2 \int_0^1 (1-t)^{3}\,\rho(t) \diff t.
	\end{align*}
	Let us now consider some explicit choices for the weight function.
	\begin{itemize}
		\item In the simplest case $\rho(x) \equiv 1$ one recovers 
		\[
		\mu=\frac{(2k)!(2k+1)!}{(k!)^2} \quad \text{and} \quad u(x)=x^k(1-x)^k, 
		\]
		which is in agreement with \cite{Hindov}.
		
		\item For $\rho(x) = \frac{1}{b-a}\chi_{[a,b]}(x)$ (with $[a,b] \subseteq [0,1]$) and $k=1$ one computes 
		\[
		\mu=\frac{12}{4(2a+b)-3(a+b)^2}.
		\]
		Now, one can choose a family of weight functions shrinking to a point, leading to the case
		\item $\rho(x) = \delta(x-a)$, for which 
		\[
		\mu=\frac{(2k-1)((k-1)!)^2}{(a(1-a))^{2k-1}}
		\]
		is the optimal constant for pointwise estimates, in agreement with \cite{Kalyabin1}. The corresponding minimizer is given by
		\[
		u(x)=
		\left\{
		\begin{aligned}
			&(1-a)^kx^kH(1-x,1-a),\quad  x\in [0,a], \\[2mm]
			&a^k(1-x)^kH(x,a),\quad  x\in [a,1],
		\end{aligned}
		\right.
		\]
		where $H(x,a)=\sum_{n=0}^{k-1}x^n\sum_{m=0}^n\binom{2k-1}{m}\binom{k-1+n-m}{n-m}a^{k-1-m}$.  
		\item We can also take $\rho = 1/x^k$ on $(0,1)$ to obtain sharp Hardy type estimates
		\[
		\int_{0}^{1}\frac{|u(x)|}{x^k}\diff x\leqslant C(k)\left(\int_{0}^{1}|u^{(k)}(x)|^2\diff x\right)^\frac12.
		\]
		For $k = 1$ one has $C(k) = 1$ and the minimizer is given by $u(x) = -x \ln x$.
	\end{itemize}

	\section{Proof of Theorem \ref{thm:1}}
	
	The existence of a minimizer $u \in W^{2k,2}(0,1)$ is provided by a standard compactness argument, and we omit it here. Computing the variation of \eqref{inf} one finds that $u$ solves the following eigenvalue problem subject to a constraint and Dirichlet boundary condition,
	\begin{equation*}
		\left\{
		\begin{aligned}
			(-1)^\ord u^{(2\ord)}(x) &=\mu \sgn(u(x)) \rho(x),  &&x\in [0,1], \\[2mm]
			\int_{0}^{1}|u(x)| \rho(x)  \diff x &= 1, \\
			u^{(j)}(0) = u^{(j)}(1) &= 0,  &&j = 0,\dots,k-1.
		\end{aligned}
		\right.
	\end{equation*}
	
	Our goal is to show that $u$ does not change its sign inside $(0,1)$. Our argument here is simply based on the following maximum principle.
	
	\begin{lemma} \label{lemma}
		Let $w$ be a solution to 
		\[
		\left\{
		\begin{aligned}
			(-1)^\ord w^{(2\ord)}(x) &=f,  &&x\in [0,1], \\[2mm]
			w^{(j)}(0) = w^{(j)}(1) &= 0,  &&j = 0,\dots,k-1.
		\end{aligned}
		\right.
		\]
		with some non-negative function $f \in L^1[0,1]$. Then $w > 0$ on $(0,1)$, unless it is zero identically.
	\end{lemma}
	
	Assuming the lemma is proved, we proceed as follows. Let $v$ be a positive eigenfunction with eigenvalue $\mu = \lambda$. Note that such a solution is well-defined and unique. Indeed, equation $(-1)^\ord v^{(2\ord)}(x) = \lambda \rho$ for a given $\lambda$ has a unique solution that satisfies the Dirichlet boundary conditions. Then the value of $\lambda$ is uniquely defined by the normalization constraint. The corresponding solution $v$ has sign by Lemma \ref{lemma} and can be chosen as positive.
	
	Then we define
	\[
	w = u - \frac{\mu}{\lambda} v
	\]
	and apply Lemma \ref{lemma} to get $u < \frac{\mu}{\lambda} v$ in $(-1,1)$, unless $w = v$ identically and we are done. Just as before, we put
	\[
	w = u + \frac{\mu}{\lambda} v
	\]
	and apply Lemma \ref{lemma} one more time to conclude $u > - \frac{\mu}{\lambda} v$. All together, that shows that
	\[
	|u| < \frac{\mu}{\lambda} v.
	\]
	Now we multiply this inequality by $\rho$, integrate it over $[0,1]$, and obtain $\mu > \lambda$. This leads to a contradiction, which shows that either $u = \frac{\mu}{\lambda} v$ or $u = -\frac{\mu}{\lambda} v$. In both cases $u$ has sign and then $\lambda = \mu$.
	
	\subsection{Proof of Lemma \ref{lemma}}    
	
	Let us show that $w$ is strictly negative on $(0,1)$, provided that it is not identically zero. Assume to the contrary that $w(x_0) = 0$ for some $x_0 \in (0,1)$. This implies that there are two zeros of $w'$ inside $(0,1)$. Arguing this way, we can show that $w^{(k)}$ must have $k+1$ distinct zeros inside $(0,1)$. This, in its turn, shows that $w^{(2k-1)}$ must have at least two zeros inside $(0,1)$. Let us consider the maximal interval $I = [a,b]$ such that $w^{(2k-1)}(a) = w^{(2k-1)}(b) = 0$. Then $w^{(2k-1)}$ is identically zero in $I$, because 
	\[
	(-1)^k w^{(2k)} = f \geqslant 0
	\]
	on $[0,1]$. Thus, $w^{(2k-2)}$ must be constant in $I$, but $I$ contains at least one zero of $w^{(2k-2)}$ so that $w^{(2k-2)} = 0$ identically in $I$. Iterating this argument, we find that $w = 0$ identically in $I$. Then $w = 0$ everywhere on $[0,1]$, which completes the proof. 
	
	\bibliographystyle{plain}
	
\end{document}